\documentclass[12pt]{article}%
\usepackage{graphicx}
\usepackage[intlimits]{amsmath}
\usepackage{amsfonts}
\usepackage{amssymb}%
 \makeatletter
\newtheorem{theorem}{Theorem}

\newtheorem{corollary}[theorem]{Corollary}

\newtheorem{problem}[theorem]{Problem}

\newenvironment{proof}[1][Proof]{\noindent{\textbf {#1}  }}  {\hfill$\Box$\bigskip}

\begin{document}

\title{On the sum of $k$ largest singular values of graphs and matrices}
\author{Vladimir Nikiforov\thanks{Department of Mathematical Sciences, University of
Memphis, Memphis TN 38152, USA} \thanks{Research supported by NSF Grant
DMS-0906634.}}
\maketitle

\begin{abstract}
In the recent years, the trace norm of graphs has been extensively studied
under the name of graph energy. The trace norm is just one of the Ky Fan
$k$-norms, given by the sum of the $k$ largest singular values, which are
studied more generally in the present paper. Several relations to chromatic
number, spectral radius, spread, and to other fundamental parameters are
outlined. Some results are extended to more general matrices.\bigskip

\textbf{Keywords: }\textit{Ky Fan norms; graph energy; singular values;
Hadamard matrices.}

\end{abstract}

\section{Introduction}

In this paper we study extremal properties of Ky Fan norms of adjacency
matrices of graphs and of matrices in general. We write $\mathcal{M}_{m,n}$
for the set of complex matrices of size $m\times n,$ and $A^{\ast}$ for the
Hermitian adjoint of a matrix $A.$ Given integers $n\geq m\geq k\geq1$ and a
matrix $A\in\mathcal{M}_{m,n},$ the Ky Fan $k$-norm $\left\Vert A\right\Vert
_{F_{k}}$ of $A\ $is defined by%
\[
\left\Vert A\right\Vert _{F_{k}}=\sigma_{1}\left(  A\right)  +\cdots
+\sigma_{k}\left(  A\right)  ,
\]
where $\sigma_{1}\left(  A\right)  \geq\sigma_{2}\left(  A\right)  \geq\cdots$
are the singular values of $A,$ that is to say, the square roots of the
eigenvalues of $AA^{\ast}.$

If $G$ is a graph, we set for short $\left\Vert G\right\Vert _{F_{k}%
}=\left\Vert A\right\Vert _{F_{k}},$ where $A$ with adjacency matrix of $G.$
Since the singular values of a Hermitian matrix are the moduli of its
eigenvalues, if $G$ is a graph of order $n,$ the parameter $\left\Vert
G\right\Vert _{F_{n}}$ is the well-studied \emph{energy} of $G,$ introduced by
Gutman in \cite{Gut78}. It is somewhat surprising that in the abundant
literature on graph energy, it hasn't been noted that the energy of a graph
$G$ is just the \emph{trace} or \emph{nuclear} norm of its adjacency matrix.
We note that this norm is widely studied in matrix theory and functional
analysis; thus, it seems that graph energy is interesting precisely because
the trace norm is a fundamental matrix parameter anyway.

Below we shall show that some problems and results in spectral graph theory
are best stated in terms of the Ky Fan norms, for example, these norms are
related to energy, spread, spectral radius, and other parameters. Thus, we
suggest to study arbitrary Ky Fan norms of graphs, in addition to the energy.
In particular, the following general problem seems interesting:

\begin{problem}
Study the extrema of $\left\Vert G\right\Vert _{F_{k}},$ and their relations
to the structure of $G$.
\end{problem}

In this note we extend and improve several results along this line of
research. On the other hand, many sound results about graphs can be readily
extended to matrices, sometimes even to non-square ones. Such facts prompt
another line of investigation:

\begin{problem}
Adopting techniques from graph theory, study extremal properties of
$\left\Vert A\right\Vert _{F_{k}},$ and their relations to the structure of
$A$ when $A$ belongs to a given class of matrices.
\end{problem}

We\ give three such results in Section \ref{sec 2}, but they are just the tip
of the iceberg.

The rest of the paper is organized as follows: in Sections \ref{sec 1} we
discuss upper bounds on Ky Fan norms of graphs; in particular, we extend a
recent result of Mohar on the sum of the largest eigenvalues, and a lower
bound on energy due to Caporossi et al. In Section \ref{sec 2}, we extend some
of the results in Section \ref{sec 1} to matrices as general as possible. At
the end we outline some open problems.

\section{Main results}

For general graph theoretic and matrix notation we refer the reader to
\cite{Bol98} and \cite{HoJo88}. Given a matrix $A=\left[  a_{ij}\right]
\in\mathcal{M}_{m,n},$ we set $\left\vert A\right\vert _{\infty}=\max
_{i,j}\left\vert a_{ij}\right\vert ,$ and $\left\vert A\right\vert
_{2}=\left(  \sum_{i,j}\left\vert a_{ij}\right\vert ^{2}\right)  ^{1/2}.$ We
say that a matrix $A$ is \emph{plain} if the all one vectors $\mathbf{j}%
_{m}\in\mathbb{R}^{m}$ and $\mathbf{j}_{n}\in\mathbb{R}^{n}$ are singular
vectors to $\sigma_{1}\left(  A\right)  ,$ that is to say, $\sigma_{1}\left(
A\right)  =\left\langle \mathbf{j}_{m},A\mathbf{j}_{n}\right\rangle /\sqrt
{mn}.$ Also $J_{m,n}$ stands for the all ones matrix of size $m\times n$, and
$J_{n}$ stands for $J_{n,n}.$

Let us first note a well-known relation that we shall use further: for every
$A\in M_{m,n},$ we have
\begin{equation}
\sigma_{1}^{2}\left(  A\right)  +\cdots+\sigma_{m}^{2}\left(  A\right)
=tr\left(  AA^{\ast}\right)  =\sum_{i,j}\left\vert a_{ij}\right\vert
^{2}=\left\vert A\right\vert _{2}^{2}. \label{eq0}%
\end{equation}

Note also that the Ky Fan norms are unitarily invariant matrix norms, a
property that may be useful in some applications. We shall use a weaker
consequence of this fact, namely, if $A$ and $B$ are matrices of size $m\times
n$, then
\[
\left\Vert A+B\right\Vert _{F_{k}}\leq\left\Vert A\right\Vert _{F_{k}%
}+\left\Vert B\right\Vert _{F_{k}}%
\]
for all $k=1,\ldots,\min\left(  m,n\right)  .$

Given a graph $G$ of order $n,$ let $\mu_{1}\left(  G\right)  ,\ldots,\mu
_{n}\left(  G\right)  $ be the eigenvalues of the adjacency matrix of $G$ in
non-increasing order.

\subsection{\label{sec 1}Upper bounds on Ky Fan norms of graphs}

In this subsection we study the asymptotics of the maximal Ky Fan $k$-norms of
graphs. To approach the problem, let us define the functions $\tau_{k}\left(
n\right)  $ and $\xi_{k}\left(  n\right)  $ as
\begin{align*}
\tau_{k}\left(  n\right)   &  =\max_{v\left(  G\right)  =n}\mu_{1}\left(
G\right)  +\cdots+\mu_{k}\left(  G\right)  ,\\
\xi_{k}\left(  n\right)   &  =\max_{v\left(  G\right)  =n}\left\Vert
G\right\Vert _{F_{k}}.
\end{align*}

For large $n$ the function $\tau_{k}\left(  n\right)  $ is pretty stable, as
implied by Theorem 1 in \cite{Nik06}:

\begin{theorem}
For every fixed positive integer $k,$ the limit $\tau_{k}=\lim
\limits_{n\rightarrow\infty}\tau_{k}\left(  n\right)  /n$ exists.
\end{theorem}

Following the approach of \cite{Nik06},\ one can prove an analogous assertion
for singular values:

\begin{theorem}
For every fixed positive integer $k,$ the limit $\xi_{k}=\lim
\limits_{n\rightarrow\infty}\xi_{k}\left(  n\right)  /n$ exists.
\end{theorem}

Admittedly, finding $\tau_{k}\left(  n\right)  $ and $\xi_{k}\left(  n\right)
$ is not easy for any $k\geq2,$ and even finding the limits $\tau_{k}$ and
$\xi_{k}$ is challenging. Indeed, even the simplest case $\xi_{2}$ is not
known yet, despite intensive research. Here is the story: in \cite{GHK01},
Gregory, Hershkowitz and Kirkland asked what is the maximal value of the
spread of a graph of order $n,$ that is to say, what is
\[
\max_{v\left(  G\right)  =n}\mu_{1}\left(  G\right)  -\mu_{n}\left(  G\right)
.
\]
This problem is still open, even asymptotically; we will not solve it here,
but we will show that it is equivalent to finding $\xi_{2}\left(  n\right)  .$
Indeed, since the singular values of a real symmetric matrix are the moduli of
its eigenvalues, we have
\[
\left\Vert G\right\Vert _{F_{2}}=\max\left\{  \left\vert \mu_{1}\left(
G\right)  \right\vert +\left\vert \mu_{2}\left(  G\right)  \right\vert
,\left\vert \mu_{1}\left(  G\right)  \right\vert +\left\vert \mu_{n}\left(
G\right)  \right\vert \right\}  ,
\]
and from \cite{EMNA08} it is known that every graph $G$ of order $n$
satisfies
\[
\left\vert \mu_{1}\left(  G\right)  \right\vert +\left\vert \mu_{2}\left(
G\right)  \right\vert \leq\left(  1/2+\sqrt{5/12}\right)  n<1.146n.
\]
On the other hand, in \cite{GHK01}, for every $n\geq2,$ a graph $G$ of order
$n$ is constructed such that
\[
\mu_{1}\left(  G\right)  -\mu_{n}\left(  G\right)  \geq\left(  2n-1\right)
/\sqrt{3};
\]
hence, $\left\vert \mu_{1}\left(  G\right)  \right\vert +\left\vert \mu
_{n}\left(  G\right)  \right\vert >1.154n$ for sufficiently large $n.$
Therefore, for large $n$ we see that
\[
\max_{v\left(  G\right)  =n}\mu_{1}\left(  G\right)  -\mu_{n}\left(  G\right)
=\max_{v\left(  G\right)  =n}\left\Vert G\right\Vert _{F_{2}}=\xi_{2}\left(
n\right)  ,
\]
as claimed.

\subsubsection{The asymptotics of $\tau_{k}\left(  n\right)  $ and $\xi
_{k}\left(  n\right)  $}

Surprisingly, finding $\tau_{k}$ and $\xi_{k}$ can be more successful for
large $k$. Indeed, Mohar \cite{Moh09} proved the following bounds%
\begin{equation}
\frac{1}{2}\left(  \frac{1}{2}+\sqrt{k}-o\left(  k^{-2/5}\right)  \right)
\leq\tau_{k}\leq\frac{1}{2}\left(  1+\sqrt{k}\right)  . \label{Mohbo}%
\end{equation}

Here we first strengthen the upper bound in (\ref{Mohbo}) and extend it to
arbitrary $\left(  0,1\right)  $-matrices, as follows:

\begin{theorem}
\label{tNik}Let $n\geq$ $m\geq k\geq1$ be integers. If $A$ is a $\left(
0,1\right)  $-matrix of size $m\times n,$ then%
\begin{equation}
\left\Vert A\right\Vert _{F_{k}}\leq\frac{1}{2}\left(  1+\sqrt{k}\right)
\sqrt{mn}. \label{in1}%
\end{equation}
Equality holds in (\ref{in1}) if and only if the matrix $J_{m,n}-2A$ is plain
and has exactly $k$ nonzero singular values which are equal. In particular, if
$m=n=k,$ equality holds if and only if $J_{m}-2A$ is a plain Hadamard matrix.
\end{theorem}

\begin{proof}
Note that $J_{m,n}-2A$ is a $\left(  -1,1\right)  $-matrix and so, in view of
the AM-QM inequality and (\ref{eq0}), we find that%
\begin{align*}
\sum_{i=1}^{k}\sigma_{i}\left(  J_{m,n}-2A\right)   &  \leq\sqrt{k\sum
_{i=1}^{k}\sigma_{i}^{2}\left(  J_{m,n}-2A\right)  }\leq\sqrt{k\sum_{i=1}%
^{n}\sigma_{i}^{2}\left(  J_{m,n}-2A\right)  }\\
&  =\sqrt{k\left\vert J_{m,n}-2A\right\vert _{2}^{2}}=\sqrt{kmn}.
\end{align*}
Therefore, using the fact that $\left\Vert \cdot\right\Vert _{F_{k}}$ is a
norm, we see that
\[
2\left\Vert A\right\Vert _{F_{k}}=\left\Vert 2A\right\Vert _{F_{k}}%
\leq\left\Vert 2A-J_{m,n}\right\Vert _{F_{k}}+\left\Vert J_{m,n}\right\Vert
_{F_{k}}\leq\sqrt{kmn}+\sqrt{mn},
\]
completing the proof of (\ref{in1}).

Clearly equality in (\ref{in1}) holds if and only if
\[
\sum_{i=1}^{k}\sigma_{i}\left(  J_{m,n}-2A\right)  =\sqrt{k\sum_{i=1}%
^{k}\sigma_{i}^{2}\left(  J_{m,n}-2A\right)  }=\sqrt{k\left\vert
J_{m,n}-2A\right\vert _{2}^{2}},
\]
and this can happen only if $\sigma_{1}\left(  J_{m,n}-2A\right)
=\cdots=\sigma_{k}\left(  J_{m,n}-2A\right)  $ and all other singular values
of $J_{m,n}-2A$ are $0$. Set $\sigma=\sigma_{1}\left(  J_{m,n}-2A\right)  .$
By the singular value interlacing theorem we have
\begin{equation}
\sigma_{i+1}\left(  2A\right)  \leq\sigma_{i}\left(  2A-J_{m,n}\right)
+\sigma_{2}\left(  J_{m,n}\right)  =\sigma\label{int}%
\end{equation}
for all $i=1,\ldots,k-1.$ Also, since $\sigma_{1}\left(  \cdot\right)  $ is a
norm, we have
\begin{equation}
\sigma_{1}\left(  2A\right)  \leq\sigma_{1}\left(  2A-J_{m,n}\right)
+\sigma_{1}\left(  J_{m,n}\right)  =\sigma+\sqrt{mn}. \label{sig}%
\end{equation}
Hence, equality in (\ref{in1}) implies equality in (\ref{int}) and
(\ref{sig}). Now let $\mathbf{x}\in\mathbb{R}^{m},$ $\mathbf{y}\in
\mathbb{R}^{n}$ be unit singular vectors to $\sigma_{1}\left(  2A\right)  .$
Since $A$ is nonnegative, we can choose $\mathbf{x}$ and $\mathbf{y}$
nonnegative too. We have
\begin{align*}
\sigma_{1}\left(  2A\right)   &  =\left\vert \left\langle \mathbf{x,}%
2A\mathbf{y}\right\rangle \right\vert \leq\left\vert \left\langle
\mathbf{x,}\left(  2A-J_{m,n}\right)  \mathbf{y}\right\rangle \right\vert
+\left\vert \left\langle \mathbf{x,}J_{m,n}\mathbf{y}\right\rangle \right\vert
\\
&  \leq\sigma_{1}\left(  2A-J_{m,n}\right)  +\sigma_{1}\left(  J_{m,n}\right)
=\sigma_{1}\left(  2A\right)  ,
\end{align*}
and so $\mathbf{x}$ and $\mathbf{y}$ are also singular vectors to $\sigma
_{1}\left(  J_{m,n}\right)  $ and to $\sigma_{1}\left(  2A-J_{m,n}\right)  .$
Since the singular vectors of $\sigma_{1}\left(  J_{m,n}\right)  $ are scalar
multiples of $\mathbf{j}_{m}$ and $\mathbf{j}_{n},$ we see that
\begin{equation}
\mathbf{x}=\left(  m^{-1/2},\ldots,m^{-1/2}\right)  ,\text{ }\mathbf{y}%
=\left(  n^{-1/2},\ldots,\alpha n^{-1/2}\right)  . \label{dxy}%
\end{equation}
Therefore, $2A-J_{m,n}$ is plain.

Now let $J_{m,n}-2A$ be a plain matrix, let $\sigma_{1}\left(  J_{m,n}%
-2A\right)  =\cdots=\sigma_{k}\left(  J_{m,n}-2A\right)  $ and all other
singular values of $J_{m,n}-2A$ be $0$. We have to prove that there is
equality in (\ref{in1}). We see immediately that $\left\Vert 2A-J_{m,n}%
\right\Vert _{F_{k}}=\sqrt{kmn},$ so to finish the proof we need to prove
that
\[
\left\Vert 2A\right\Vert _{F_{k}}=\left\Vert 2A-J_{m,n}\right\Vert _{F_{k}%
}+\left\Vert J_{m,n}\right\Vert _{F_{k}}.
\]
Let $\mathbf{x}$ and $\mathbf{y}$ be defined by (\ref{dxy}). We have
\[
\sigma_{1}\left(  2A-J_{m,n}\right)  =\left\langle \mathbf{x},\left(
2A-J_{m,n}\right)  \mathbf{y}\right\rangle \text{ \ \ and \ \ }\sigma
_{1}\left(  J_{m,n}\right)  =\left\langle \mathbf{x},J_{m,n}\mathbf{y}%
\right\rangle ,
\]
and therefore,
\[
\sigma_{1}\left(  2A\right)  \geq\left\langle \mathbf{x},2A\mathbf{y}%
\right\rangle =\left\langle \mathbf{x},\left(  2A-J_{m,n}\right)
\mathbf{y}\right\rangle +\left\langle \mathbf{x},J_{m,n}\mathbf{y}%
\right\rangle \geq\sigma_{1}\left(  2A\right)  ,
\]
implying that $2A$ is plain and that $\sigma_{1}\left(  2A\right)  =\sigma
_{1}\left(  2A-J_{m,n}\right)  +\sigma_{1}\left(  J_{m,n}\right)  $. It
follows that for every $i=2,\ldots,m,$ $\sigma_{i}\left(  2A\right)
=\sigma_{i}\left(  2A-J_{m,n}\right)  $ and there exist $\mathbf{x}%
_{i},\mathbf{y}_{i}$ such that
\[
\left\langle \mathbf{x}_{i},\left(  2A-J_{m,n}\right)  \mathbf{y}%
_{i}\right\rangle =\left\langle \mathbf{x}_{i},2A\mathbf{y}_{i}\right\rangle
.
\]
Indeed let us check this assertion for $i=2.$ By the singular value
decomposition theorem, we have
\begin{align*}
\sigma_{2}\left(  2A\right)   &  =\sigma_{1}\left(  2A-\sigma_{1}\left(
2A\right)  \mathbf{x}\otimes\mathbf{y}\right)  =\sigma_{1}\left(
2A-\sigma_{1}\left(  2A\right)  \frac{1}{\sqrt{mn}}J_{m,n}\right) \\
&  =\sigma_{1}\left(  2A-\left(  \sigma_{1}\left(  2A-J_{m,n}\right)
+\sigma_{1}\left(  J_{m,n}\right)  \right)  \frac{1}{\sqrt{mn}}J_{m,n}\right)
\\
&  =\sigma_{1}\left(  2A-J_{m,n}-\sigma_{1}\left(  2A-J_{m,n}\right)
\mathbf{x}\otimes\mathbf{y}\right)  =\sigma_{2}\left(  2A-J_{m,n}\right)  .
\end{align*}
Setting $\mathbf{x}_{2},\mathbf{y}_{2}$ to be singular vectors to $\sigma
_{1}\left(  2A-\sigma_{1}\left(  2A\right)  \mathbf{x}\otimes\mathbf{y}%
\right)  ,$ the assertion follows.

This completes the proof of Theorem \ref{tNik}.
\end{proof}

For graphs we get the following consequence.

\begin{corollary}
\label{tMoh}If $n$ and $k$ are integers such that $n\geq k\geq1,$ then
\[
\xi_{k}\left(  n\right)  \leq\frac{1}{2}\left(  1+\sqrt{k}\right)  n.
\]

\end{corollary}

Since for every Hermitian matrix $A$ of size $n\geq k$ we have%
\[
\mu_{1}\left(  A\right)  +\cdots+\mu_{k}\left(  A\right)  \leq\left\Vert
A\right\Vert _{F_{k}},
\]
it follows that $\tau_{k}\left(  n\right)  \leq\xi_{k}\left(  n\right)  ,$ and
so Corollary \ref{tMoh} implies the upper bound in (\ref{Mohbo}). The
advantage of Theorem \ref{tNik} is that this bound is extended to non-square
matrices, where eigenvalues are not applicable at all; moreover, along the
same lines, Theorem \ref{mo3} below gives an analogous statement for arbitrary
nonnegative matrices.

Let us also point out that Theorem \ref{tNik} neatly proves and generalizes
the result of Koolen and Moulton \cite{KoMo01}: \emph{If }$G$\emph{ is a graph
of order }$k,$\emph{ then}%
\begin{equation}
\left\Vert G\right\Vert _{F_{k}}\leq\frac{1}{2}\left(  1+\sqrt{k}\right)  k.
\label{KoMobo}%
\end{equation}

\subsubsection{Lower bounds on $\xi_{k}\left(  n\right)  $}

Theorem \ref{tNik} can be complemented by a lower bound which stems from
Mohar's lower bound on $\tau_{k}$.

\begin{theorem}
\label{tNiklo}Given $\varepsilon>0,$ for sufficiently large $k,$ $m$ and $n,$
there exists a $\left(  0,1\right)  $-matrix $A$ of size $m\times n$ such
that
\[
\left\Vert A\right\Vert _{F_{k}}\geq\frac{1}{2}\left(  \frac{1}{2}+\sqrt
{k}-\varepsilon k^{-2/5}\right)  \sqrt{mn}.
\]

\end{theorem}

\begin{proof}
Mohar in \cite{Moh09} uses a class of strongly regular graphs, to show that
for every $\varepsilon>0,$ if $k$ is sufficiently large, then
\[
\tau_{k}>\frac{1}{2}\left(  \frac{1}{2}+\sqrt{k}-\varepsilon k^{-2/5}\right)
.
\]
Fix $\varepsilon>0,$ and choose $k$ and $l$ such that
\[
\tau_{k}\left(  l\right)  >\frac{1}{2}\left(  \frac{1}{2}+\sqrt{k}%
-\frac{\varepsilon}{2}k^{-2/5}\right)  l,
\]
that is to say, there exists a graph $G$ of order $l$ such that
\[
\left\Vert G\right\Vert _{F_{k}}=\xi_{k}\left(  l\right)  \geq\tau_{k}\left(
l\right)  >\frac{1}{2}\left(  \frac{1}{2}+\sqrt{k}-\frac{\varepsilon}%
{2}k^{-2/5}\right)  l.
\]
Suppose that $m$ and $n$ are sufficiently large and set $p=\left\lfloor
m/l\right\rfloor ,$ $q=\left\lfloor n/l\right\rfloor .$ Writing $\otimes$ for
the Kronecker product and $B$ for the adjacency matrix of $G$, let $A=B\otimes
J_{p,q}.$ Since $A$ is a $pl\times ql$ matrix satisfying
\[
\left\Vert A\right\Vert _{F_{k}}=\left\Vert B\otimes J_{p,q}\right\Vert
_{F_{k}}=\left\Vert B\right\Vert _{F_{k}}\sqrt{pq},
\]
for sufficiently large $m$ and $n$ we find that
\begin{align*}
\frac{1}{\sqrt{mn}}\left\Vert A\right\Vert _{F_{k}}  &  =\frac{1}{\sqrt{mn}%
}\left\Vert B\right\Vert _{F_{k}}\sqrt{pq}>\frac{1}{2\sqrt{mn}}\left(
\frac{1}{2}+\sqrt{k}-\frac{\varepsilon}{2}k^{-2/5}\right)  l\sqrt{pq}\\
&  =\frac{1}{2}\left(  \frac{1}{2}+\sqrt{k}-\frac{\varepsilon}{2}%
k^{-2/5}\right)  \sqrt{\left(  1-\frac{l}{m}\right)  \left(  1-\frac{l}%
{n}\right)  }\\
&  >\frac{1}{2}\left(  \frac{1}{2}+\sqrt{k}-\frac{\varepsilon}{2}%
k^{-2/5}\right)  \sqrt{1-\frac{l}{m}-\frac{l}{n}}\\
&  >\frac{1}{2}\left(  1+\sqrt{k}-\varepsilon k^{-2/5}\right)  .
\end{align*}
To complete the proof it is enough to make $A$ an $m\times n$ matrix by adding
$n-ql$ zero columns and $m-pl$ zero rows$.$ Clearly these additions do not
affect the singular values of $A.$
\end{proof}

The proof of Theorem \ref{tNiklo}, in fact, implies the following corollary
for graphs.

\begin{corollary}
\label{tMohlo}Given $\varepsilon>0,$ for all sufficiently large $k$ and $n,$
\[
\xi_{k}\left(  n\right)  \geq\frac{1}{2}\left(  \frac{1}{2}+\sqrt
{k}-\varepsilon k^{-2/5}\right)  n.
\]

\end{corollary}

The bounds in Theorem \ref{tNiklo} and Corollary \ref{tMohlo}, although more
flexible that the lower bound in (\ref{Mohbo}), are essentially equivalent to
it. However, there are infinitely many cases when $\xi_{k}\left(  n\right)  $
behaves better and attains the upper bound of Corollary \ref{tMoh}. Indeed,
let $k$ be such that there exists a graph $G$ of order $k$ for which the
Koolen and Moulton bound (\ref{KoMobo}) is an equality. Then, for every
$n\geq1,$ blowing-up $G$ by a coefficient $n,$ and calculating the Ky Fan
$k$-norm of the resulting graph, we obtain%
\[
\xi_{k}\left(  kn\right)  =\frac{1}{2}\left(  1+\sqrt{k}\right)  kn.
\]
In particular, it is known that for $k=4$ the complete graph $K_{4}$ satisfies%
\[
\left\Vert K_{4}\right\Vert _{F_{4}}=6=\frac{1}{2}\left(  1+\sqrt{4}\right)
4,
\]
and so, $\xi_{4}\left(  4n\right)  =6n$ and\ $\xi_{4}=3/2.$ Other known $k$ of
this type are $k=4m^{4}$ for all integer $m>1,$ (see \cite{HaXi09}), but there
are also others, like $k=16$ and $k=36.$

Clearly, following the blow-up idea, for infinitely many triples $k,m,n$ one
can construct matrices attaining equality in (\ref{in1})$.$

\subsection{The Ky Fan norms and chromatic number}

In this section we first obtain an inequality which extends a result of
Caporossi, Cvetkovi\'{c}, Gutman and Hansen \cite{CCGH99}: \emph{If }$G$
\emph{is a graph of order }$n,$\emph{ then}
\begin{equation}
\left\Vert G\right\Vert _{F_{n}}\geq2\mu_{1}\left(  G\right)  , \label{Cap}%
\end{equation}
\emph{where equality holds if an only if }$G$\emph{ is a complete multipartite
graph with possibly some isolated vertices.}

Inequality (\ref{Cap}) can be improved if we know the chromatic number $\chi$
of $G.$ Indeed, letting $n$ be the order of $G,$ recall that Hoffman's
inequality \cite{Hof70} gives
\[
\left\vert \mu_{n}\left(  G\right)  \right\vert +\cdots+\left\vert \mu
_{n-\chi+2}\left(  G\right)  \right\vert \geq\mu_{1}\left(  G\right)  .
\]
From here we easily obtain the following theorem, improving (\ref{Cap}):

\begin{theorem}
\label{tHof}If $G$ is a graph with chromatic number $\chi\geq2,$ then%
\begin{equation}
\left\Vert G\right\Vert _{F_{\chi}}\geq2\mu_{1}\left(  G\right)  .
\label{Hofin}%
\end{equation}

\end{theorem}

Note that if $G$ is a complete $\chi$-partite graph with possibly some
isolated vertices, then equality holds in (\ref{Hofin}). However, there are
many other cases of equality some of which are rather complicated and their
complete description seems difficult.

In contrast to Theorem \ref{tHof}, for bipartite graphs we have
\[
\left\Vert G\right\Vert _{F_{2}}=2\mu_{1}\left(  G\right)  \leq n.
\]
It is not clear how the above inequality can be generalized to $r$-partite
graphs for $r>2.$

Moreover, since every triangle-free graph of order $n$ satisfies $\mu
_{1}\left(  G\right)  \leq\sqrt{m},$ see \cite{Nos70} or \cite{Nik02}, we
obtain the following

\begin{theorem}
\label{tTfree}If $G$ is a triangle-free graph of order $n,$ then%
\[
\left\Vert G\right\Vert _{F_{2}}\leq2\sqrt{m}.
\]
Equality holds if and only if $G$ is a complete bipartite graph with possibly
some isolated vertices.
\end{theorem}

It is not clear how Theorem \ref{tTfree} can be generalized to $K_{r}$-free
graphs for $r>3.$

\subsection{\label{sec 2}The maximal Ky Fan norms of matrices}

Here our aim is to generalize Theorem \ref{tNik} to arbitrary matrices.

First, applying the AM-QM inequality to the sum of the largest $k$ singular
values and using (\ref{eq0}), we obtain the following theorem.

\begin{theorem}
\label{mo1}Let $n\geq m\geq2,$ and $m\geq k\geq1.$ If $A\in M_{m,n},$ then
\[
\left\Vert A\right\Vert _{F_{k}}\leq\sqrt{k}\left\vert A\right\vert _{2}.
\]
Equality holds if and only if $A$ has exactly $k$ nonzero singular values and
they are equal.
\end{theorem}

It seems difficult to characterize effectively the class of all matrices with
the above property. Here we give a construction which suggests the great
diversity of this class. Let $q\geq k$ and let $B$ be a $k\times q$ matrix
whose rows are pairwise orthogonal vectors of length $c$. Since $BB^{\ast
}=c^{2}I_{k}$, we see that all $k$ singular values of $B$ are equal to $c$.
Setting $A=B\otimes J_{r,s},$ we see that $A$ has exactly $k$ nonzero singular
values, which are equal, and so $\left\Vert A\right\Vert _{F_{k}}=\sqrt
{k}\left\vert A\right\vert _{2}$.

If we use that $\left\vert A\right\vert _{2}\leq\sqrt{mn}\left\vert
A\right\vert _{\infty},$ obtaining the following theorem:

\begin{theorem}
\label{mo2}Let $n\geq m\geq2,$ and $m\geq k\geq1.$ If $A\in M_{m,n},$ then
\begin{equation}
\left\Vert A\right\Vert _{F_{k}}\leq\sqrt{kmn}\left\vert A\right\vert
_{\infty}. \label{Hadin}%
\end{equation}
Equality is possible if and only if all entries of $A$ has the same absolute
value, $A$ has exactly $k$ nonzero singular values and they are equal.
\end{theorem}

We gave this easy theorem with the sole purpose to discuss the class of
matrices for which equality holds in (\ref{Hadin}). To this end, let $n\geq
m\geq2$ be integers and write $\mathsf{Had}_{m,n}$ for the class of $m\times
n$ matrices whose entries have the same absolute value and whose rows are
pairwise orthogonal. Clearly if $A\in\mathsf{Had}_{m,n},$ then all row vectors
of $A$ have the same length. Note also that the class $\mathsf{Had}_{n,n}$ is
just the set of all scalar multiples of Hadamard matrices.

Here is a general construction of matrices for which equality holds in
(\ref{Hadin}): Let $q\geq k,$ let $B\in\mathsf{Had}_{k,q},$ and set
$A=B\otimes J_{r,s}.$ Since $A$ has exactly $k$ nonzero singular values and
they are equal, and since the entries of $A$ are equal in absolute value, we
have $\left\Vert A\right\Vert _{F_{k}}=\sqrt{k}\left\vert A\right\vert
_{2}=\sqrt{kmn}\left\vert A\right\vert _{\infty},$ thus equality holds for $A$
in (\ref{Hadin}).

Although the upper bounds given in Theorems \ref{mo1} and \ref{mo2} are as
good as one can get, for non-negative matrices there is a slight improvement.

\begin{theorem}
\label{mo3}Let $n\geq m\geq2,$ and $m\geq k\geq1.$ If $A\in M_{m,n}$ is a
nonnegative matrix, then
\[
\left\Vert A\right\Vert _{F_{k}}\leq\frac{1}{2}\left(  1+\sqrt{k}\right)
\sqrt{mn}\left\vert A\right\vert _{\infty}.
\]
Equality holds in (\ref{in1}) if and only if $A$ is a scalar multiple of a
$\left(  0,1\right)  $-matrix, the matrix $\left\vert A\right\vert _{\infty
}\left(  J_{m,n}-2A\right)  $ is plain and has exactly $k$ nonzero singular
values which are equal.
\end{theorem}

The proof of Theorem \ref{mo3} essentially repeats the proof of Theorem
\ref{tNik} and we will omit it.\bigskip

\textbf{Some open problems}\bigskip

(1) Find $\xi_{k}\left(  n\right)  $ and $\tau_{k}\left(  n\right)  ;$

(2) Find the best approximation of $\xi_{n}\left(  n\right)  $ for all $n;$

(3) Find the extrema of $\left\Vert G\right\Vert _{F_{k}}$ when $G$ belongs to
a given monotone or hereditary property.

(4) Same problem for matrices.

(5) Characterize the graphs for which equality holds in (\ref{Hofin});

(6) Characterize effectively the matrices giving equality in Theorems
\ref{mo1} and \ref{mo2}.

\end{document}